\documentclass[12pt]{amsart}
\usepackage{verbatim}
\usepackage{graphicx}

\title{Constant mean curvature surfaces of 
any positive genus}
\author{S-P. Kobayashi}
\address{S-P. Kobayashi \& W. Rossman, 
Department of Mathematics,
Kobe University, Rokko Kobe 657-8501, Japan} 
\email{kobayasi@math.kobe-u.ac.jp}  
\email{wayne@math.kobe-u.ac.jp}  

\author{M. Kilian}
\address{M. Kilian, 
Department of Mathematical Sciences,
University of Bath, Bath, BA2 7AY, United Kingdom} 
\email{masmk@maths.bath.ac.uk}  

\author{W. Rossman}

\author{N. Schmitt}
\address{N. Schmitt, Tech. Univ. Berlin,
Fakult\"{a}t II, Mathematik und Naturwissenshaften,
Stra{\ss}e des 17.~Juni 136, 10623 Berlin, Germany} 
\email{nick@gang.umass.edu}

\numberwithin{equation}{section}

\newcommand{\abs}[1]{{\lvert#1\rvert}}

\newcommand{\half}{\tfrac{1}{2}}

\newcommand{\ol}[1]{\overline{#1}}

\DeclareMathOperator{\re}{Re}

\DeclareMathOperator{\suchthat}{|}

\DeclareMathOperator{\tr}{tr}

\theoremstyle{plain}
\newtheorem{theorem}{Theorem}[section]
\newtheorem{corollary}[theorem]{Corollary}
\newtheorem{lemma}[theorem]{Lemma}

\newtheorem*{corollary*}{Corollary}
\newtheorem*{lemma*}{Lemma}
\newtheorem*{proposition*}{Proposition}
\newtheorem*{theorem*}{Theorem}
\theoremstyle{definition}

\newtheorem{definition}[theorem]{Definition}
\newtheorem{example}{Example}

\newtheorem{remark}{Remark}[section]
\newtheorem*{algorithm*}{Algorithm}
\newtheorem*{application*}{Application}
\newtheorem*{assertion*}{Assertion}
\newtheorem*{assumption*}{Assumption}
\newtheorem*{axiom*}{Axiom}
\newtheorem*{claim*}{Claim}
\newtheorem*{conjecture*}{Conjecture}
\newtheorem*{definition*}{Definition}
\newtheorem*{example*}{Example}
\newtheorem*{notation*}{Notation}
\newtheorem*{note*}{Note}
\newtheorem*{observation*}{Observation}
\newtheorem*{question*}{Question}
\newtheorem*{remark*}{Remark}
\theoremstyle{plain}

\newcommand{\bbC}{\mathbb{C}}

\newcommand{\bbR}{\mathbb{R}}
\newcommand{\bbS}{\mathbb{S}}

\newcommand{\calD}{\mathcal{D}}

\newcommand{\calI}{\mathcal{I}}

\newcommand{\calP}{\mathcal{P}}

\newcommand{\calT}{\mathcal{T}}

\newcommand{\SL}{\mathrm{SL}_{\mbox{\tiny{$2$}}}(\bbC)}

\newcommand{\SU}{\mathrm{SU}_{\mbox{\tiny{$2$}}}}
\newcommand{\su}{\mathfrak{su}_{\mbox{\tiny{$2$}}}}

\newcommand{\GL}{\mathrm{GL}_{\mbox{\tiny{$2$}}}(\bbC)}

\newcommand{\gl}{\mathfrak{gl}_{\mbox{\tiny{$2$}}}(\bbC)}

\newcommand{\Id}{\mathrm{Id}}

\newcommand{\be}{\begin{equation}} 
\newcommand{\ee}{\end{equation}}

\begin{document}

\thanks{Kobayashi supported by DFG grant DO 776/1. \\
\indent Kilian supported by EPSRC grant 
GR/S28655/01. \\ 
\indent Rossman supported by Japan Monbusho 
grant (B)(1)--15340023. \\
\indent 2000 {\it Mathematics Subject Classification.} 53A10.}


\begin{abstract}
\footnotesize
We show the existence of several new families 
of non-compact constant mean curvature surfaces: 
(i) singly-punctured surfaces of 
arbitrary genus $g \geq 1$, 
(ii) doubly-punctured tori, and  
(iii) doubly periodic surfaces 
with Delaunay ends.
\end{abstract}

\maketitle
\section{Introduction}

Since Wente's discovery \cite{Wen:hopf} of tori 
there has been revived interest in the study of 
nonminimal constant mean curvature ({\sc{cmc}}) 
surfaces. The investigation has been informed by 
analytical methods such as work by Kapouleas  
\cite{Kap1} and Korevar, Kusner and Solomon 
\cite{KorKS} as well as by techniques from 
integrable systems such as in Pinkall and Sterling 
\cite{PinS} and Bobenko \cite{Bob:tor}. 

In the late 1990's Dorfmeister, Pedit and Wu 
\cite{DorPW} formulated a Weierstra{\ss} 
type representation for {\sc{cmc}} surfaces, 
which involves solving 
a holomorphic complex linear $2 \times 2$ system of 
ordinary differential equations 
with values in a loop group, and subsequently 
factorising the solution into two factors, 
one of which turns out to be a moving frame of 
the Gau{\ss} map, from which the surface can be 
constructed. Unfortunately, the factorisation 
is not explicit and much qualitative information 
about the solution is obscured in the process. 
In representing surfaces with non-trivial 
topology the main challenge is to keep track 
of the monodromy representation of 
the moving frame. 

Any {\sc{cmc}} surface comes in an isometric 
$\mathbb{S}^1$ family, but the periods of a 
non-simply connected surface are generally 
only closed for one 
specific value $\lambda_0 \in \mathbb{S}^1$ 
of the \emph{spectral parameter}. 
The Weierstra{\ss} data for a {\sc{cmc}} surface 
consists of a Riemann surface $\Sigma$, a point 
$z_0$ on the universal cover $\widetilde{\Sigma}$, 
a \emph{holomorphic potential} $\xi$ on $\Sigma$ 
and an initial condition $\Phi_0$. 
Solving the initial value problem 
\begin{equation} \label{eq:IVP}
  d\Phi = \Phi\,\xi,\,\Phi(z_0) = \Phi_0
\end{equation}
yields a solution $\Phi$ and corresponding 
monodromy representation both of which depend 
on the spectral 
parameter, since generally both $\xi$ and 
$\Phi_0$ do. If the monodromy $M$ satisfies for 
all deck transformations 
the following three conditions 
\begin{align}
 &\left.M\right|_{\mathbb{S}^1} \in \SU, 
	\label{eq:unitary} \\
 &\left.M\right|_{\lambda_0} = \pm \Id, 
	\label{eq:per1}\\
 &\left.d_\lambda M\right|_{\lambda_0} = 0, 
	\label{eq:per2} 
\end{align}
then the resulting associated family factors 
through the fundamental group at $\lambda_0$ and 
we thus have a {\sc{cmc}} immersion 
$f:\Sigma \to \mathbb{R}^3$. In equation \eqref{eq:per2} 
and throughout this work we denote by $d_{-}$ 
the derivative with respect to the subscript, 
which we omit in the case of the exterior 
derivative on the Riemann surface, as in 
\eqref{eq:IVP}. 
Condition \eqref{eq:per1} kills the rotational 
periods while \eqref{eq:per2} takes care of the 
translational periods, and both 
can be ensured by properties on $\xi$. 
The condition \eqref{eq:unitary} 
is harder to satisfy and makes use of varying the 
initial condition $\Phi_0$. These three conditions 
have been used in a number of papers, starting with 
the work of Dorfmeister and Haak \cite{DorH:cyl} 
and by several of the authors while investigating 
{\sc{cmc}} immersions of the $n$-punctured Riemann 
sphere, the so called $n$-Noids \cite{KilKRS}, 
\cite{KilSS} and \cite{Sch:tri}. 
Another approach to studying embedded {\sc{cmc}} 
$3$-Noids can be found in the work of 
Gro{\ss}e-Brauckmann, Kusner and Sullivan 
\cite{GroKS:Tri}.  

The purpose of this paper is to show the existence 
of new {\sc{cmc}} surfaces by exhibiting 
Weierstra{\ss} data 
$(\,\Sigma,\,\xi,\,\Phi_0,\,z_0\,)$ 
that fulfill the above requirements 
\eqref{eq:unitary}--\eqref{eq:per2} and to initiate 
the study of higher genus surfaces via loop group 
techniques. Briefly summarising the contents of 
this paper, after providing 
some general sufficient conditions on 
Weierstra{\ss} 
data to satisfy the condition \eqref{eq:unitary}, 
\eqref{eq:per1} and \eqref{eq:per2} we apply these 
results to prove existence of new examples of 
{\sc{cmc}} surfaces
\begin{enumerate} 
\item of any positive genus and a single end,   
\item of genus $1$ with two ends, 
\item which are doubly-periodic with infinitely 
	many ends asymptotic to Delaunay ends.  
\end{enumerate}
Although the last mentioned surfaces (iii) are 
immersions of genus zero domains with 
infinitely many punctures, 
they have natural quotient surfaces with 
positive genus.    

\section{Preliminary results}

We denote an annular neighbourhood of 
the unit circle $\mathbb S^1$ for some real 
$r \in (0,\,1]$ by 
$A_r = \left\{ \lambda \in \mathbb C : 
  r \leq |\lambda | \leq 1/r \right\}$.
It is common abuse to call a map $M: A_r \to \SL$ 
\emph{unitary} if $\left.M\right|_{\mathbb S^1} \in \SU$.
\begin{definition}
We shall call a map $M: A_r \to \SL$ 
\emph{unitarisable} on $A_s$ for some 
  $s \in [r,\,1]$ if there exists 
  a map $h : A_s \to \GL$ for some 
  $s \in [r,\,1]$ such that 
  $h\,M\,h^{-1}: A_s \to \SL$ is unitary.
\end{definition}
We use the following notation for diagonal and 
off--diagonal $2 \times 2$ matrices:
\begin{equation*}
\mathrm{diag}[u,\, v] = 
\bigl( \begin{smallmatrix} 
u&0\\0&v \end{smallmatrix} \bigr),\,
\mathrm{off}[u,\, v] = 
\bigl( \begin{smallmatrix} 
0&u\\v&0 \end{smallmatrix} \bigr). 
\end{equation*}
In preparation for Theorems \ref{thm:puncture1} and 
\ref{thm:doublypunctured} we first provide some 
technical results. The next lemma gives conditions 
on a matrix which ensure that after unitarisation, 
it satisfies the closing conditions at $\lambda_0 = 
e^{i\,x_0}$.
\begin{lemma}
  \label{thm:unitrace}
  Let $J\subset\mathbb R$ be an open interval, and 
  let $M:J\to\SL,\,U:J \to \SU$ be smooth maps
  with $\tr M=\tr U$. If 
  $M(x_0) = \pm \Id$ and $d_xM (x_0)$ 
  is nilpotent for $x_0 \in J$, then 
  $U(x_0)= \pm \Id$ and $d_x U(x_0)=0$.  
\end{lemma}
\begin{proof}
  Since $\tr U(x_0)=\pm 2$ and $U(x_0)\in\SU$, 
  we have $U(x_0)=\pm \Id$.
  Let $\tau=\tfrac{1}{2}\tr M = \tfrac{1}{2}\tr U$.  We 
  differentiate the Cayley-Hamilton equations 
  $M^2-2 \tau M = U^2-2 \tau U = -\Id$ twice 
  and evaluate at $x_0$ to get
  $\pm{d_x M(x_0)}^2 = d_x^2\tau (x_0)\Id = 
  \pm{d_x U(x_0)}^2$. 
  So $d_x M(x_0)$ nilpotent implies $d_x U(x_0)$ 
  is also nilpotent.
  Since $U \in \SU$, we have $d_x U \in \su$, 
  and so nilpotency implies 
  $d_x U(x_0)=0$.
\end{proof}
Lemma \ref{thm:dtrace} computes the derivatives 
of a solution to a linear ODE with respect to 
a parameter. From this the series
expansion of the trace of the monodromy 
with respect to the parameter can be computed, 
and hence the trace
can be estimated in a small interval, 
as we will see in Lemma \ref{thm:dtrace3}.  
See~\cite{DorH:cyl} for a related theorem. 
\begin{lemma} \label{thm:dtrace}
  Let $\Sigma$ be a simply connected 
  Riemann surface with coordinate $z$.  
  Let $A(z,\,x): \Sigma \times \mathbb R \to \gl$ 
  be analytic in $z$ and smooth in $x$ and 
  $B(x): \mathbb R \to \gl$ be smooth. 
  Let $X(z,\,x): \Sigma 
  \times \mathbb R \to \gl$ 
  be the solution of the initial value problem
  \begin{equation} \label{eq:phi0}
    \left(d_z X \right)(z,\,x) =
    X(z,\,x) A(z,\,x)\, , \quad X(z_0,\,x) = B(x).
  \end{equation}
  Let $X_k(z): \Sigma \to \gl$, for integers 
  $k \geq 0$, be 
  the solutions to the sequence of initial 
  value problems 
\begin{equation*}
  \left( d_z X_k \right)(z) = \sum_{i,j 
  \geq 0,i+j=k}
  \frac{k!}{i!(k-i)!}
  X_i(z) A_j(z) \, , \quad  X_k(z_0) = B_k,
\end{equation*}
  where 
  $A_j(z) := \left(d^j_x A \right)(z,\,x_0)$ and 
  $B_k := \left(d_x^k B\right)(x_0)$. Then
  \begin{equation}\label{eq:phiKpart}
  \left( d_x^k X \right)(z,\,x_0) = X_k(z).
  \end{equation}
\end{lemma}
\begin{proof}
  Differentiate~\eqref{eq:phi0} 
  repeatedly with respect to $x$.
\end{proof}
In the following, let $\Sigma$ be a connected 
Riemann surface with universal cover 
$\widetilde{\Sigma}$ and $\Delta$ its group of 
deck transformations. We denote the 
holomorphic $1$-forms on $\Sigma$ by 
$\Omega'(\Sigma,\,\mathbb C)$. 

In the next Lemma we show that a certain class of 
potentials always ensures the closing conditions 
\eqref{eq:per1} and \eqref{eq:per2}. 
Such potentials will be used in later examples 
(Theorems \ref{thm:puncture1} and 
\ref{thm:doublypunctured}) to show 
the existence of new {\sc{cmc}} surfaces.
\begin{lemma} \label{thm:dtrace2}
  Let $f,\,g \in \Omega'(\Sigma,\,\mathbb C)$ and 
  $t = \lambda^{-1}(\lambda-1)^2$ and 
  \begin{equation} \label{eq:firstA}
    A = \begin{pmatrix}
    0 & f \, t \\ g & 0
    \end{pmatrix} \; . 
  \end{equation}
  Let $w_0 \in \widetilde{\Sigma}$ and $X$ 
  be the solution to the initial value problem
  \begin{equation} \label{eq:phi1}
    d X  = X\, A\,, \quad X(w_0,\,t) = \Id \; .
  \end{equation}
  Let $\gamma \in \Delta$ and 
  $M(t):= X(\gamma(w_0),\,t)$.  
  Suppose that
  \begin{equation} \label{eq:a2}
    \int_{w_0}^{\gamma(w_0)} g = 0.
  \end{equation}
  Then $\tilde{M}(1)=\Id$ and 
  $d_\lambda \tilde{M} (1)=0$, 
  where $\tilde{M}(\lambda) = M(t)$. 
\end{lemma}

\begin{proof}
  Note that 
  $X(w,\,0) = \Id + 
  \mathrm{off}[\,0,\, \int_{w_0}^w g\,]$. 
  Hence $X(\gamma(w_0),\,0)=\Id$, and so 
  $\tilde{M}(1)=\Id$.  
  Then with $z_0=w_0$, $e^{ix}=\lambda$, 
  $e^{ix_0}=\lambda_0=1$, $B(x)=\Id$, 
  $A$ as in \eqref{eq:firstA} and 
  $z$ fixed to $\gamma(w_0)$ in 
  \eqref{eq:phiKpart}, it follows 
  that $A_1(z)$ is identically zero, and 
  Lemma~\ref{thm:dtrace} implies that 
  $(d_\lambda \tilde{M})(1)=0$.
\end{proof}

The next lemma will be used in the proofs of 
Theorems \ref{thm:puncture1} and 
\ref{thm:doublypunctured} to show that 
certain monodromy groups can be 
unitarised.  

\begin{lemma}
  \label{thm:dtrace3}
  Take the same notations and conditions as in 
  Lemma~\ref{thm:dtrace2},
  with $t$ replaced by $ct$ for 
  some constant $c \in \mathbb R \setminus \{ 0 \}$.  
  Suppose that $\tau(ct) = \tfrac{1}{2}\tr M(ct)$
  is real for all $t \in [-4,0]$ and that 
  \begin{equation} \label{eq:a4}
    I_0 I_2+I_1^2 < 0 \; , 
  \end{equation} 
  where 
$$
  I_0 = \int_{w_0}^{\gamma(w_0)}\hspace{-3mm} f\,,
  	\quad 
  I_1 = \int_{w_0}^{\gamma(w_0)}\hspace{-2mm}
	\left(f\int g\right) \quad  
  \mbox{and} \quad 
  I_2 = \int_{w_0}^{\gamma(w_0)}\hspace{-2mm} 
  \left( g  \int \left( f \int g \right) \right).
$$
  Then for $\tilde{\tau}(\lambda,\,c)= 
  \tau(c\,t)=
  \tau(c\,\lambda^{-1} (\lambda-1)^2)$ 
  there exists a $c_0>0$ such that for 
  all $|c| \in (0,\, c_0)$ we have  
\begin{equation*}
  \abs{\tilde{\tau}(\lambda,\,c)}<1 
  \text{ for all } \lambda \in 
  \mathbb{S}^1\setminus\{1\} \;. 
\end{equation*}
\end{lemma}
\begin{proof}
  Let $X$ be the solution of 
  equation \eqref{eq:phi1} with $f$ 
  replaced by $cf$, 
  and write $X_{ij} = X_{ij}(w,\,t,\,c)$ for the 
  entries of $X$. Defining 
  $\hat{X}_{ij} = \tfrac{d}{d(ct)}X_{ij}
  (\gamma(w_0),\, 0,\, c)$, 
  Lemma~\ref{thm:dtrace} implies 
  \begin{equation*} 
	\begin{split}
    \hat{X}_{11} &= \int_{w_0}^{\gamma(w_0)} 
               \hspace{-2mm}\left(g \int f\right),\,
    \hat{X}_{12} = \int_{w_0}^{\gamma(w_0)}  
               \hspace{-6mm}f \\ 
    \hat{X}_{21} &= \int_{w_0}^{\gamma(w_0)} 
               \hspace{-2mm}\left(g \int f \left( 
	\int g \right) \right) \mbox{  and  }
    \hat{X}_{22} = \int_{w_0}^{\gamma(w_0)}
               \hspace{-2mm}\left( f \int g \right). 
	\end{split}
  \end{equation*}
  Note that the $\hat{X}_{ij}$ are independent 
  of $c$. Considering the 
  first two derivatives of $\det X = 1$ 
  with respect to $ct$ and evaluating 
  at $t=0$ and $w=\gamma(w_0)$ yields 
  \begin{equation*}
  \left( d_{ct} \tau\right)(0,\,c) = 0 \; , \;\;\; 
  \left( d^2_{ct} \tau\right)(0,\,c) =
  \hat{X}_{12} \hat{X}_{21} - 
	\hat{X}_{11} \hat{X}_{22} = 
	\hat{X}_{12} \hat{X}_{21} 
  +\hat{X}_{22}^2 \; . 
  \end{equation*}
  Note that the first of these two 
  equations implies 
  $\hat{X}_{11} = - \hat{X}_{22}$, 
  which is used in the second of 
  these two equations.  
  Equation \eqref{eq:a4} implies that 
  $\hat{X}_{12} \hat{X}_{21} 
  +\hat{X}_{22}^2 < 0$, and so 
  the second derivative with respect to 
  $ct$ of $\tau$ is negative, and $\tau$ 
  attains a maximum of $1$ at $t=0$.
  Thus there exists a $k_0>0$ such that 
  $\abs{ct}\in(0,\,k_0]$ implies 
  $\abs{\tau(c t)} \in [0,\,1)$.
  Let $c_0=k_0/4$. Then for all $c$ 
  such that $|c| \in (0,c_0)$ and for all 
  $t \in [-4,0]$, we have $|\tau(c t)| \in [0,1)$ 
  and the lemma follows.  
\end{proof}

Applying Lemma \ref{thm:dtrace}, 
similarly to the proofs of 
Lemmas \ref{thm:dtrace2} and 
\ref{thm:dtrace3}, we obtain
\begin{corollary}\label{cor:prelim}
With notations and conditions as in 
Lemmas \ref{thm:dtrace2} and \ref{thm:dtrace3}, 
we have $d^2_\lambda \tilde{M} (1) = 
2 ( \mathrm{diag}[ 
  \,-I_1,\,I_1\,] + \mathrm{off}[\,I_0,\,I_2\,])$.
\end{corollary}


\section{Singly-punctured {\sc{cmc}} surfaces 
of arbitrary genus}
\begin{figure}
  \centering
  \includegraphics[scale=1]
	{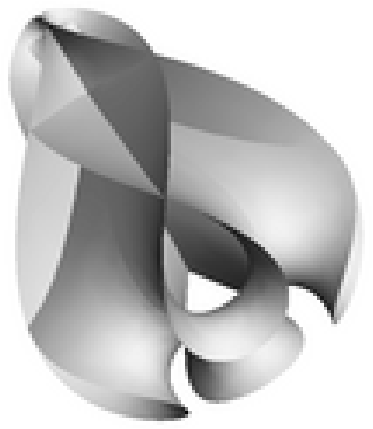}
  \includegraphics[scale=1]
	{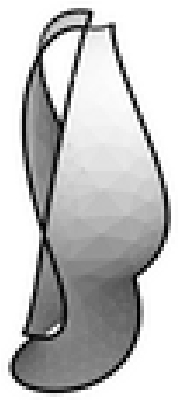}
  \includegraphics[scale=.85]
	{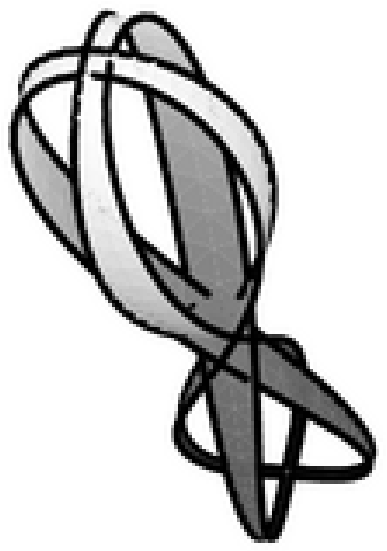}
  \includegraphics[scale=1]
	{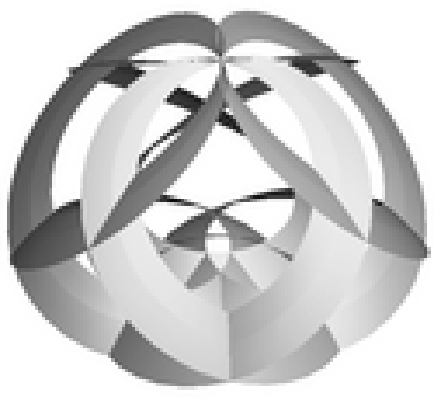}
  \caption{The three figures on the left are 
  parts of a {\sc{cmc}} singly-punctured torus, 
  as in Theorem \ref{thm:puncture1} with $n=2$. 
  The left-most image shows the torus with a 
  neighborhood of the end removed.
  The surface has $90^\circ$ rotation 
  symmetry and reflection symmetry.  
  The second image shows one-fourth of 
  the surface, which extends to the 
  full surface (again with a neighborhood 
  of the end removed) by these 
  symmetries. The third image is a skeletal 
  portion of the surface.  
  The Hopf differential 
  has a pole of order $6$ at the end. 
  The right-most figure is a 
  {\sc{cmc}} doubly-punctured torus, 
  as in Theorem~\ref{thm:doublypunctured}, and 
  the image here shows a skeleton of this 
  torus (with a doubly-punctured disk 
  containing the ends removed).}
  \label{fig:genus1}
\end{figure}

We construct a family of {\sc{cmc}} 
immersions of a singly-punctured genus 
$g$ Riemann surface into $\bbR^3$
with umbilics, for any positive $g$. 
The closing problem is solved by 
imposing symmetries so that the 
monodromy group can be shown to be unitarisable.
\begin{theorem} \label{thm:puncture1}
  Let $n\ge 2$ be an even integer.  
  Let $\Sigma$ be the singly-punctured
  hyperelliptic genus $n/2$ Riemann surface 
  defined by $\Sigma = \{(z,\,w)\in\bbC^2 
  \suchthat w^2=z(1-z^n)\}$. Let 
  \begin{equation*}
  \xi = 
  \begin{pmatrix}
    0 & c\lambda^{-1}(\lambda-1)^2 w^{-1}dz \\ 
  d(z^{n-1}w) & 0
  \end{pmatrix},\quad c \in \bbR^\ast.
  \end{equation*}
  Then for $c$ sufficiently close to zero, 
  $\xi$ induces a conformal {\sc{cmc}} 
  immersion $\Sigma\to\bbR^3$ with 
  order $2n$ dihedral symmetry.  
\end{theorem}
\begin{proof}
  Choose a basepoint 
  $(z_0,\,w_0) \in \widetilde{\Sigma}$ 
  in the fibre of $(0,\,0) \in \Sigma$ and let
  $\Phi$ be the solution to the initial 
  value problem \eqref{eq:IVP} with 
  $\Phi_0 = \Phi(z_0,\,w_0) = \Id$.
  We must show that there exists a 
  unitariser for the monodromy of $\Phi$ 
  and verify the closing conditions.  

  With $\alpha=\exp(\pi i/n)$,
  for $k\in\{0,\dots,n-1\}$
  let $\gamma_k:[0,\,1]\to\Sigma$ be the 
  curve from $(0,\,0)$ to
  $(\alpha^{2 k},\,0)$ and back to $(0,\,0)$ 
  along the straight
  line in the $z$-plane from $0$ to 
  $\alpha^{2 k}$, defined by 
  \begin{equation*} \gamma_k(s) = \left\{
	\begin{tabular}{ll}
	$\left( 2 s \alpha^{2k},
  	(-\alpha)^{k} \left| 
  	\sqrt{2 s (1-2^n s^n)}\right| \right)$, 
	&$0 \leq s \leq 1/2$ \\
  	$\left( 2 (1-s) 
	\alpha^{2k},- (-\alpha)^{k} 
        \left| \sqrt{2 (1-s) (1-2^n (1-s)^n)} 
      	\right| \right)$ 
	&$1/2 \leq s \leq 1$. \end{tabular}
	\right.
  \end{equation*}
  Let $\hat{\gamma}_k$ be the lifted curves 
  originating at $(z_0,\,w_0)$ and 
  $M_k(\lambda) := 
  \Phi(\hat{\gamma}_k(1),\,\lambda)$.
  Then $M_0,\dots,M_{n-1}$ generate the 
  monodromy group of $\Phi$, since 
  $\Sigma$ has only one puncture at $(z,\, w)=
  (\infty,\, \infty)$, and the monodromy 
  about the puncture is 
  $N=\prod_{k=0}^{n-1} M_k$. 
  Lemma~\ref{thm:dtrace2} then implies 
  \begin{equation}\label{eq:closingbeforedressing}
    M_k(1)=\Id,\quad
    d_\lambda M_k(1)=0 \; , 
    \quad k\in\{0,\dots,n-1\} \; .  
  \end{equation}
  Hence we also have $N(1)=\Id$, $d_\lambda N(1)=0$. 
  It remains to show that there exists a unitariser 
  for the monodromy group. For this, we compute 
  the monodromy group's symmetries. We define 
  the following maps on $\Sigma$: 
  \begin{equation} \label{eq:symmetry}
    \sigma(z,\,w)=(\alpha^2 z,\,\alpha w)\; ,\, 
    \rho(z,\,w)=(z,\,-w)\; \mbox{ and }
    \theta(z,\,w)=(\ol{z},\,\ol{w}) \; . 
  \end{equation}
  Then for $g_\sigma = 
	\mathrm{diag}[\,\sqrt{\alpha}^{-1},\,
	\sqrt{\alpha}\,]$ and 
  $g_\rho = \mathrm{diag}[\,-i,\,i\,]$ we have 
\begin{equation*}
  \sigma^\ast\xi = 
  g_\sigma^{-1}\,\xi\,g_\sigma,\, 
  \rho^\ast\xi = g_\rho^{-1}\,\xi\,g_\rho 
  \mbox{ and }
  \ol{\theta^\ast\xi(1/\bar{\lambda})} = 
  \xi(\lambda),
\end{equation*}
  where expressions like $\sigma^\ast\xi$ denote 
  $\sigma^\ast\xi((z,\,w),\,\lambda) = 
  \xi(\sigma(z,\,w),\,\lambda)$.
  Since $(0,\,0)$ is 
  a fixed point of $\sigma,\,\rho$ and $\theta$, 
  we define the lifts $\hat\sigma,\,\hat\rho,\,
  \hat\theta$ that map 
  $(0,\,0)$ to $(z_0,\,w_0)$. Then using 
  $\Phi(z_0,\,w_0)=\Id$, we obtain 
\begin{equation*} 
    \hat\sigma^\ast\Phi = 
	g_\sigma^{-1}\Phi \,g_\sigma \; , \;\;\; 
    \hat\rho^\ast\Phi = 
	g_\rho^{-1}\Phi \,g_\rho \; , \;\;\; 
    \ol{\hat\theta^\ast\Phi(1/\bar{\lambda})} 
	= \Phi(\lambda) \; . 
\end{equation*}  
  Hence the monodromy group has the following 
  symmetries:
  \begin{equation} \label{eqn:monodromysymmetries}
    \begin{split}
    &M_k^{(-1)^k} = 
	g_\sigma^{-k}M_0 \,g_\sigma^k\;, 
 	\quad k\in\{0,\dots,n-1\} \; , \\
    &M_0^{-1} = g_\rho^{-1}M_0 \,g_\rho \; , \\
    & \ol{M_0} = M_0 \;\;\;\; 
	\text{for all} \;\, \lambda \in \bbS^1 \; . 
    \end{split}
  \end{equation}
  Note that the third of these symmetries also 
  follows from the facts that 
  $\lambda^{-1} (\lambda-1)^2 \in \bbR$ 
  for all $\lambda \in \bbS^1$ and 
  $\gamma_0(s) \in \bbR^2$ for all $s \in [0,1]$. 
  Denoting the entries of $M_0$ by $M_{ij}$,  
  the third symmetry in 
  \eqref{eqn:monodromysymmetries} implies that 
  the $M_{ij}$ are all real, and the second 
  symmetry in 
  \eqref{eqn:monodromysymmetries} implies that 
  $M_{11}=M_{22}$, for all $\lambda \in \bbS^1$.  
  In particular, $\tr (M_0)$ is 
  real for all $\lambda \in \bbS^1$. 
  The integrals $I_j$ in equation~\eqref{eq:a4} 
  are then 
  \begin{equation*}
  I_0 = 2 c \int_{\hat{\gamma}} |w|^{-1}dz\;,\quad 
  I_1 = 0 \; , \quad \mbox{ and }
  I_2 = \frac{2 c}{n} \int_{\hat{\gamma}} 
  z^n d(z^{n-1} |w|) \; , 
  \end{equation*}
  for the curve $\hat{\gamma}(s) = 
  (s,|\sqrt{s(1-s^n)}|)\in\Sigma$, $s \in [0,1]$.  
  Using the formula, valid for 
  $\re n>0,\ \re r>0$ and $\re s>0$, 
  \begin{equation*}
  \int_{\hat{\gamma}} z^{r-1}(1-z^n)^{s-1}\,dz=
  \frac{\Gamma(\tfrac{r}{n})\,\Gamma(s)}
  {n\,\Gamma(\tfrac{r}{n}+s)} \; , 
  \end{equation*}
  where $\Gamma(\ell)=\int_0^\infty y^{\ell-1} 
  e^{-y} dy$ is the Euler gamma function, 
  we get 
  \begin{equation*}
  I_0 I_2+I_1^2 = -\frac{8 \pi c^2 (2n-1)
  \cot(\tfrac{\pi}{2n})}{(n-1)(3n-1)(5n-1)} < 0.
  \end{equation*}
  Hence by Lemma~\ref{thm:dtrace3},
  with $\tau(\lambda,\, c)=\half\tr (M_0(\lambda))$, 
  there exists a $c_0>0$ such that
  for all $c$ satisfying $\abs{c}\in(0,\,c_0)$,
  \begin{equation}\label{thetaubound}
  \abs{\tau(\lambda,\,c)}<1 \text{ for all } 
  \lambda\in\bbS^1\setminus\{ 1 \} \; , 
  \end{equation} and $\tau(1,\, c)=1$.  
  Then $\tau=M_{11}=M_{22} \in \bbR$ has 
  modulus at most $1$ for all $\lambda \in 
  \bbS^1$.  
  Thus $-M_{12}M_{21}=1-M_{11}^2 \geq 0$ on 
  $\bbS^1$, so 
\begin{equation*} 
  v:=-\frac{M_{21}}{M_{12}} 
  \geq 0 \text{ on } \bbS^1 \; . 
\end{equation*}
  Furthermore, $v$ is finite and strictly 
  positive on 
  $\bbS^1 \setminus \{ 1 \}$, by 
  \eqref{thetaubound}.  
  
  Let us now consider the behavior 
  of $v$ at $\lambda=1$. 
  By \eqref{eq:closingbeforedressing}, 
  we know that 
\begin{equation*} 
  M_{12}|_{\lambda=1}=M_{21}|_{\lambda=1}= 
  d_\lambda M_{12}|_{\lambda=1}=
  d_\lambda M_{21}|_{\lambda=1}=0 \; . 
\end{equation*}
  Applying Corollary \ref{cor:prelim}, 
  we have $d_\lambda^2 M_{12}|_{\lambda=1}=2I_0$ 
  and $d_\lambda^2 M_{21}|_{\lambda=1}=2I_2$. 
  Since $I_1=0$ and $I_0I_2+I_1^2 < 0$, we conclude 
  that $I_0$ and $I_2$ are both nonzero, 
  so $M_{12}$ and $M_{21}$ 
  both have zeroes of order exactly 
  two at $\lambda = 1$. Hence $v$ is nonzero and 
  finite at $\lambda = 1$.  Thus $v$ is a strictly 
  positive finite function on all of $\bbS^1$, 
  and therefore $\sqrt[4]{v}$ can be 
  globally and smoothly defined on $\bbS^1$. 
  Then by the first symmetry of 
  \eqref{eqn:monodromysymmetries}, 
  the diagonal unitariser 
  is given by 
\begin{equation*} 
  h = \begin{pmatrix} \sqrt[4]{v} & 0 \\ 
  0 & (\sqrt[4]{v})^{-1} \end{pmatrix}, 
\end{equation*}
  which simultaneously unitarizes 
  $M_0,\dots,M_{n-1}$ on $\bbS^1$, i.e. 
  $h M_j h^{-1} \in \SU$ for all 
  $\lambda \in \bbS^1$. 
  Therefore the monodromy group of 
  $h \Phi$ is unitarized on all of $\bbS^1$.  

  By Equation \eqref{eq:closingbeforedressing} 
  and Lemma~\ref{thm:unitrace}, 
  the monodromy group $h M_j h^{-1}$ still 
  satisfies the closing 
  conditions \eqref{eq:per1} and \eqref{eq:per2} 
  at $\lambda=1$.
  Hence the resulting {\sc{cmc}} immersion is 
  well-defined on $\Sigma$.  

  Since the coefficient $c w^{-1} dz$ of 
  the $\lambda^{-1}$ term of the 
  upper-right entry of the potential 
  $\xi$ has no zeros or poles on the 
  singly-punctured Riemann
  surface $\Sigma$, the {\sc{cmc}} 
  immersion is unbranched, 
  see \cite{DorH:cyl}, Theorem 3.1.  

  We now consider the symmetries of the 
  {\sc{cmc}} immersion resulting from $h \Phi$.  
  Since $(0,\,0)$ is fixed by the map 
  $\sigma$ and $h$ is 
  independent of $(z,w)$ and $[h,\,g_\sigma]=0$ 
  and also $\hat{\sigma}^\ast\Phi = 
  g_\sigma^{-1}\Phi g_\sigma$, 
  we have that $(\hat{\sigma}^k)^\ast (h \Phi) = 
  g_\sigma^{-k} (h \Phi) 
  g_\sigma^k$, where $\hat{\sigma}^k$ is the 
  composition of $\hat{\sigma}$ with itself 
  $k$ times.  

  Let $h \Phi = F B$ be the Iwasawa 
  decomposition with respect to $\bbS^1$ 
  (Theorem 8.1.1 \cite{PreS}), pointwise 
  on $\widetilde{\Sigma}$. Since 
  $g_\sigma^{-k} F g_\sigma^k$ is unitary and 
  $g_\sigma^{-k} B g_\sigma^k$ positive, 
  the unitary part of 
  $(\hat{\sigma}^k)^\ast (h_1 \Phi)$ is 
  $g_\sigma^{-k} F g_\sigma^k$. 
  The symmetry $(\hat{\sigma}^k)^\ast F = 
    g_\sigma^{-k} F g_\sigma^k$ 
  descends to the immersion via the Sym-Bobenko 
  formula \cite{Bob:tor}, see also \cite{KilKRS} 
  section 4, and results 
  in a rotation of angle $k\pi/n$ about 
  an axis independent of $k$.  
  Hence the surface has an order 
  $2n$ rotational symmetry.  

  To show dihedral symmetry, we now need only 
  show that the surface has at least 
  one reflective symmetry across a 
  plane parallel to the common 
  axis of the rotational symmetries.  
  We will show that the map 
  $\theta(z,w)=(\bar z,\bar w)$ 
  is such a reflective symmetry, by 
  showing that the immersion generated by $F$ 
  via the Sym-Bobenko formula \cite{Bob:tor}, 
  and denoted by $f$, satisfies
  \begin{equation*} 
    \theta^\ast f = - \bar f \; . 
  \end{equation*}  
  Because $\xi|_{\lambda} = \xi|_{\lambda^{-1}}$, 
  we have $\Phi(z,w,\lambda) = 
  \Phi(z,w,\lambda^{-1})$ and consequently   
  \begin{equation*}  
    \overline{\Phi(\bar z,\bar w,\bar \lambda)} = 
    \overline{\Phi(\bar z,\bar w,\bar \lambda^{-1})} 
    = \overline{\hat{\theta}^\ast 
	\Phi(\bar \lambda^{-1})} 
    = \Phi(z,w,\lambda).
  \end{equation*}
  This further implies that 
  $\overline{\Phi(\gamma_0(s),\bar \lambda)} = 
  \Phi(\gamma_0(s),\lambda)$ and so 
  $\overline{M(\bar \lambda)} = M(\lambda)$, 
  and in turn 
  $\overline{h(\bar \lambda)} = h(\lambda)$, 
  since $\theta^\ast \gamma_0(s) = 
  \gamma_0(s)$. Thus 
  \begin{equation*}
    \overline{h(\bar \lambda) 
    \,\Phi(\bar z,\bar w,\bar \lambda)} = 
    h(\lambda) \Phi(z,w,\lambda)
  \end{equation*}
  and consequently  
  $\overline{F(\bar z,\bar w,\bar \lambda)} \, 
    \overline{B(\bar z,\bar w,\bar \lambda)} = 
    F(z,w,\lambda) \, B(z,w,\lambda)$. 
  Uniqueness of the Iwasawa decomposition 
  yields $\overline{F(\bar z,\bar w,\bar \lambda)} = 
  F(z,w,\lambda)$ and implies 
  $\theta^\ast f = - \bar f$.  
\end{proof}

\begin{remark}
  Note that the end of any surface in Theorem 
  \ref{thm:puncture1} is not 
  asymptotically Delaunay, because 
  the order of the Hopf differential 
  there is strictly less than $-2$.
  This is also implied by \cite{KorKS}, 
  since Delaunay ends have non-zero weight, 
  but the balancing formula implies that the single 
  end of any surface in Theorem 
  \ref{thm:puncture1} must have zero weight.  
\end{remark}


\section{{\sc{cmc}} immersions of a 
doubly-punctured torus}

In this section, we construct immersions 
of a doubly-punctured genus $1$ 
Riemann surface into $\bbR^3$ with umbilics.  
\begin{theorem}\label{thm:doublypunctured}
  Let $\calT = \{ [z] \in \bbC / \Gamma \, | \, 
  z \in \bbC \}$ 
  be the square torus, 
  where $\Gamma$ is the $2$-dimensional lattice 
  generated by $2 \omega_1 \in \bbR^+$ and 
  $2 \omega_2 = 2 i \omega_1$.  
  Let $\omega_3 = \omega_1+\omega_2$. 
  On the twice-punctured torus $\Sigma = 
  \calT\setminus\{[\omega_3/2],\,[-\omega_3/2] \}$, 
  let $\xi$ be the 
  potential
  \begin{equation*}
    \xi =
    \begin{pmatrix}
    0 & c\lambda^{-1}(\lambda-1)^2 \\
    \wp''''(z+\omega_3/2)+\wp''''(z-\omega_3/2) & 0
  \end{pmatrix}dz \; , \;\;\; c \in \bbR^* , 
  \end{equation*} 
  where $\wp$ is the Weierstrass $\wp$-function 
  with respect to $\calT$ satisfying 
  $(\wp')^2=4\wp(\wp^2-1)$ and $\, '$ 
  denotes the derivative 
  with respect to $z$.  
  Then for $c$ sufficiently close to zero, 
  $\xi$ induces a conformal 
  {\sc{cmc}} immersion $\Sigma \to \bbR^3$ 
  with order $4$ dihedral symmetry.  
\end{theorem}

\begin{remark}\label{rem:aboutWP}
Note that $\wp''''= 120 \wp^3 - 72 \wp$. 
Then, since $\wp(-z)=\wp(z)$ and 
$\wp(iz)=-\wp(z)$, it follows that also 
$\wp''''(-z)=\wp''''(z)$ and 
$\wp''''(iz)=-\wp''''(z)$.  
These properties will be used 
in the following proof.  
One other particular property that 
we will need is, defining 
\begin{align*} 
  \mathcal{I} (z) &= 
     \wp'''(z+\omega_3/2)+\wp'''(z-\omega_3/2)-
     \wp'''(\omega_3/2)-\wp'''(-\omega_3/2) \\
     &= \wp'''(z+\omega_3/2)+\wp'''(z-\omega_3/2), 
\end{align*}
that the integral 
$\int_0^{2 \omega_1} (\calI (z))^2 dz > 0$ 
along the real axis from $0$ to 
$2 \omega_1$ is positive.  
This integral is real because 
of the relations $(\calI (\omega_1 \pm \bar{z}))^2
=\overline{(\calI (z))}^2$, and then one 
can check that it is positive for any choice 
of $\omega_1 > 0$.  
\end{remark}

\begin{proof}
  Choose a basepoint $w_0 \in \widetilde{\Sigma}$ 
  in the fibre of $z_0 = 0 \in \bbC$, and let 
  $\Phi$ be the solution to the initial 
  value problem 
  $d\Phi = \Phi\xi$, $\Phi(w_0)= \Id$. 
  Let $\gamma_k=\gamma_k(s) \in \bbC$ be the 
  straight-line curve from $z_0$ to $2\omega_k$ 
  ($k\in\{1,\,2\}$) defined 
  by $\gamma_k(s) = 2 s \omega_k$ for 
  $s \in [0,1]$.  
  Let $\delta_1=\delta_1(s) \in \bbC$ 
  for $s \in [0,1]$ 
  be a curve from $z_0$ around $\omega_3/2$
  in the counterclockwise direction and back to 
  $z_0$ lying in a small 
  neighborhood of the straight line 
  from $z_0$ to $\omega_3/2$, 
  and let $\delta_2=\delta_2(s)=-\delta_1(s)$ 
  be the curve from $z_0$ around $-\omega_3/2$ 
  in the counterclockwise direction and back 
  to $z_0$ that is the reflection of 
  $\delta_1$ through the point $z_0$.  

  Let $M_k=M_k(\lambda)$ be the respective 
  global monodromies of $\Phi$ over the torus 
  along $\gamma_k$, and let 
  $A_k=A_k(\lambda)$ be the monodromies of 
  $\Phi$ about the two punctures of $\Sigma$ 
  along $\delta_k$ ($k\in\{1,\,2\}$).  
  Then $M_1,\,M_2,\,A_1,\,A_2$ generate 
  the monodromy group of $\Phi$.  
  
  This proof follows the same strategy as the 
  proof of Theorem~\ref{thm:puncture1}.  
  First, we note that all the generating 
  elements $M_1,M_2,A_1,A_2$ 
  of the monodromy group of $\Phi$ satisfy the 
  closing conditions \eqref{eq:per1} and 
  \eqref{eq:per2}. This follows from 
  Lemma~\ref{thm:dtrace2}, since the 
  lower-left entry in $\xi$ is the derivative 
  with respect to $z$ of a function that is 
  well-defined on $\Sigma$ and hence 
  Equation \eqref{eq:a2} will be satisfied. 
  Our main effort again goes into showing that 
  there exists an initial condition that 
  unitarises the monodromy group of $\Phi$. 
  To accomplish this, we first 
  compute the symmetries of the monodromy group 
  and define the following 
  transformations of $\Sigma$: 
\begin{equation*} 
    \sigma(z) = z + \omega_3 \; , \;\;\; 
    \rho(z) = iz + \omega_1 \; , \;\;\; 
    \theta(z) = \bar{z} + \omega_1 \; . 
\end{equation*}
  Then with 
  $g=\mathrm{diag}[\,1/\sqrt{i},\,\sqrt{i}\,]$, 
  the potential $\xi$ has the symmetries 
  \begin{equation*} 
    \sigma^\ast\xi = \xi \; , \;\;\; 
    \rho^\ast\xi = g^{-1}\,\xi \,g \; ,\mbox{ and } 
    \ol{\theta^\ast\xi(1/\bar{\lambda})} = \xi \; . 
  \end{equation*}
  Hence 
  $\hat{\sigma}^\ast\Phi = V_\sigma \Phi$, 
  $\hat{\rho}^\ast\Phi = V_\rho \Phi g$ and 
  $\ol{\hat{\theta}^\ast \Phi (1/\bar{\lambda})} = 
  V_\theta \Phi$
  for some $z$-independent 
  $V_\sigma,\,V_\rho$ and $V_\theta$.  
  Since $z_0=0$ is a fixed point of the two maps 
  \begin{equation*} 
  \sigma^{-1}\rho^2: z\mapsto -z \; , \quad 
  \sigma^{-1}\rho \,\theta: z\mapsto i \bar{z} 
  \end{equation*}
  (we interpret these compositions 
  as being applied in order from rightmost 
  first to leftmost last), 
  and since $\Phi(z_0)=\Id$, we have 
  $V_\rho^2 V_\sigma^{-1} = g^{-2}$ and 
  $V_\theta V_\rho V_\sigma^{-1} = g^{-1}$. 
  It follows that 
  \begin{equation}
  \label{eqn:thisguytookawhile}
  \begin{split}
    M_1^{-1} = g^{-2}M_1 g^{2} \; &, \;\;\; 
    A_2 = g^{-2}A_1 g^{2} \; , \\ 
    \ol{M_2(1/\bar{\lambda})} = 
    g^{-1} M_1(\lambda) g \; &, \;\;\; 
    \ol{A_1(1/\bar{\lambda})}^{-1} = 
    g^{-1} A_1(\lambda) g \; . 
  \end{split}
  \end{equation}
  The first and third equations in 
  \eqref{eqn:thisguytookawhile} imply that 
  also $M_2^{-1}=g^{-2}M_2 g^{2}$.  

  Because the potential is real-valued along 
  the curve $\gamma_1$ when 
  $\lambda \in \bbS^1$, we conclude that $M_1$ 
  is a real-valued matrix for all 
  $\lambda \in \bbS^1$.  This 
  fact, combined with the first equation 
  in \eqref{eqn:thisguytookawhile}, implies 
  that $M_1$ has the form 
\begin{equation*}
  M_1 = 
  \begin{pmatrix}
    a_1 & b_1 \\ c_1 & a_1
  \end{pmatrix},
\end{equation*}
  where
  $a_1=a_1(\lambda)=\ol{a_1(1/\bar{\lambda})}$, \ 
  $b_1=b_1(\lambda)=\ol{b_1(1/\bar{\lambda})}$ and 
  $c_1=c_1(\lambda)=\ol{c_1(1/\bar{\lambda})}$. 
  Furthermore, the fourth equation in 
  \eqref{eqn:thisguytookawhile} implies 
  \begin{equation*}
  A_1 =
  \begin{pmatrix}
    a_2 & b_2 \\ c_2 & d_2
  \end{pmatrix} \; , \quad \text{where} 
\end{equation*}
  \begin{equation}\label{a2b2c2d2}
    \ol{a_2(1/\bar{\lambda})}=d_2(\lambda)\, ,\; 
    \ol{b_2(1/\bar{\lambda})}=-i b_2(\lambda)\,,\; 
    \ol{c_2(1/\bar{\lambda})}=i c_2(\lambda) \, . 
  \end{equation}
  From this it is clear that 
  $\tau_1 = \tfrac{1}{2} \tr M_1$ and 
  $\tau_2 = \frac{1}{2} \tr A_1$ are 
  real for all $\lambda \in \bbS^1$.  

  We will now show 
  that $M_1$ and $A_1$ are simultaneously 
  unitarizable for small $|c|$.  
  Toward this goal, we first apply Lemma 
  \ref{thm:dtrace3} to show that for small $|c|$ 
  we have $|\tau_1(\lambda)|<1$ for 
  all $\lambda \in \bbS^1 \setminus \{ 1 \}$:  
  We take $\Sigma$ and $\xi$ as in Theorem 
  \ref{thm:doublypunctured} and take $f=dz$ and 
  $g = 
  (\wp''''(z+\omega_3/2)+\wp''''(z-\omega_3/2))dz$ 
  and the curve $\gamma=\gamma_1$.  
  Then, in Lemma \ref{thm:dtrace3}, we have 
  $I_0 = 2 \omega_1 > 0$ and $I_1=0$.  
  To compute $I_2$, integration by parts 
  yields
  \begin{equation*} 
    I_2 = - \int_{\gamma_1} \left( f 
    \left( \int g \right)^2 \right)  = 
    - \int_{\gamma_1} (\calI (z))^2 dz \; , 
  \end{equation*}
   where $\calI (z)$ is as defined in 
  Remark \ref{rem:aboutWP}. Then by Remark 
  \ref{rem:aboutWP}, we have $I_2 < 0$.  
  Thus the conditions of Lemma \ref{thm:dtrace3} 
  hold and we conclude that for all 
  $c \in \bbR$ sufficiently close to $0$, we have
  \begin{equation}\label{tau1in4point1} 
    |\tau_1(\lambda)|<1 \;\;\; \text{for all} \; 
    \lambda \in \bbS^1 \setminus \{ 1 \} \; . 
  \end{equation}
  From \eqref{tau1in4point1} and the fact 
  that $b_1c_1=a_1^2-1 \in \bbR$ 
  on $\bbS^1$, we have 
  \begin{equation}\label{b1c1negative} 
    (b_1c_1) |_{\lambda=1} = 0 \;\; \text{and} \;\;   (b_1c_1) |_{\bbS^1 \setminus 
  \{ 1 \} } < 0 \; . 
  \end{equation} 
  Thus we can define a function 
  \begin{equation*} 
     v = - \frac{c_1}{b_1} 
  \end{equation*} 
  that is finite and nonzero on 
  $\bbS^1 \setminus \{ 1 \}$.  
  Furthermore, by \eqref{b1c1negative} 
  and the fact that $b_1 \in \bbR$ 
  on $\bbS^1$, we conclude that
  \begin{equation}\label{eqn:positivityofV} 
    0 < v < \infty 
  \end{equation}
  for all $\lambda \in \bbS^1 \setminus \{ 1 \}$.  
  Similar to the arguments in 
  proving Theorem~\ref{thm:puncture1}, 
  Corollary \ref{cor:prelim} implies that 
  $b_1$ and $c_1$ both have zeroes of order 
  exactly two, hence $v$ is nonzero and finite 
  at $\lambda=1$ as well.  It follows that 
  $\sqrt[4]{v} > 0$, representing the positive 
  fourth root of $v$, is globally and smoothly 
  defined on $\bbS^1$.  We define
\begin{equation}
  h = \begin{pmatrix} \sqrt[4]{v} & 0 \\
	0 & (\sqrt[4]{v})^{-1} \end{pmatrix}.
\end{equation}
  Because $b_1c_1 \leq 0$ and $\sqrt{v} \geq 0$, 
  the conjugate $h M_1 h^{-1} \in \SU$ for 
  $|\lambda  |=1$.  
  
  The image of the path $\gamma_2\delta_1$ 
  under the map $\rho$ is 
  homotopic to the path $\gamma_1^{-1}\delta_1$.  
  Hence $A_1M_2$ and $A_1M_1^{-1}$ are 
  conjugate and so have 
  the same trace. Hence 
  \begin{equation*} 
    \tr(A_1 (M_2-M_1^{-1})) = b_2 c_1 (1+i) + 
    c_2 b_1 (1-i) = 0 \; . 
  \end{equation*}
  It follows that
  \begin{equation}\label{fordefiningV}
    c_2 b_1 + i c_1 b_2 = 0 \; . 
  \end{equation}
  In \eqref{fordefiningV}, either both $b_2$ 
  and $c_2$ are identically zero, or 
  neither of them are identically zero.  
  If $b_2$ and $c_2$ are identically 
  zero, then $A_1$ is diagonal and 
  $A_1 \in \SU$ for all $\lambda \in \bbS^1$.  
  Hence we have succeeded in 
  simultaneously unitarizing both $M_1$ and 
  $A_1$ on $\bbS^1$ 
  by conjugating by $h$.  We may then 
  proceed to the final paragraph of this proof, 
  which gives the concluding 
  argument for proving 
  Theorem~\ref{thm:doublypunctured}.  
  Therefore, without 
  loss of generality, let us assume 
  that neither $b_2$ nor $c_2$ is 
  identically zero.  

  Under the assumption that 
  $b_2$ and $c_2$ are not identically zero, 
  by \eqref{fordefiningV} we also have 
\begin{equation*} 
  v=-\frac{\ol{c_1(1/\bar{\lambda})}}{b_1(\lambda)} 
   =-\frac{\ol{c_2(1/\bar{\lambda})}}{b_2(\lambda)}\; . 
\end{equation*}
  Furthermore, \eqref{a2b2c2d2} and 
  \eqref{eqn:positivityofV} 
  then imply that $b_2=r_1 (1+i)$ and 
  $c_2 = r_2 (1-i)$ with $r_1,r_2 
  \in \bbR$ and $r_1^{-1}r_2 \leq 0$, 
  on $\bbS^1$. These facts 
  together show that also the conjugate 
  $h A_1 h^{-1} \in \SU$ for 
  $|\lambda  |=1$.  
  
  Thus we have simultaneously unitarised 
  $M_1$ and $A_1$ on $\bbS^1$.  
  Since $g \in \SU$ and commutes with $h$, 
  conjugation by $h$ also unitarizes $M_2$ 
  and $A_2$, so the full monodromy group 
  is unitarized on $\bbS^1$. Now, like in 
  the proof of Theorem \ref{thm:puncture1}, 
  using Iwasawa splitting on $\bbS^1$ and 
  noting that the 
  monodromy of $h\, \Phi$ still satisfies 
  \eqref{eq:per1} and \eqref{eq:per2}, we
  conclude that the resulting {\sc{cmc}} surface 
  given by the Sym-Bobenko
  formula is defined on $\Sigma$. Finally, 
  analogous to the arguments at the end of 
  the proof of Theorem~\ref{thm:puncture1}, 
  the order $4$ dihedral symmetry of the 
  resulting {\sc{cmc}} immersions can be shown, 
  and since the coefficient $c dz$ of the 
  $\lambda^{-1}$ term of the 
  upper-right entry of the potential 
  $\xi$ has no zeros or poles on the 
  twice-punctured Riemann surface $\Sigma$, 
  the resulting {\sc{cmc}} immersion is 
  unbranched \cite{DorH:cyl}.  
  This completes the proof.  
\end{proof}


\section{Doubly-periodic {\sc{cmc}} surfaces in 
$\bbR^3$ with ends that are asymptotically Delaunay}

In this section, we provide a third class of 
Weierstra{\ss} data for which the monodromy 
can be unitarised. The potential is of interest 
to us because, although
the monodromy can be unitarised, the monodromy 
does not satisfy \eqref{eq:per1} at 
$\lambda_0 = 1$. The relaxing of
this closing condition is what allows the 
resulting {\sc{cmc}} immersions 
to extend to doubly-periodic 
surfaces (when $n=3,4,6$ in the
theorem). 
\begin{theorem}
  \label{lem:Section5first}
  Let $n\ge 3$ be an integer, and   
  define the Riemann surface 
  $\Sigma=(\bbC \setminus \mathcal{P}) 
  \cup \{\infty\}$ 
  with 
  $\mathcal{P} = \{ z \in \bbC \, | \, z^n=1 \}$.  
  Let 
\begin{equation} \label{eq:potentiallemma}
    \xi = 
    \begin{pmatrix}  0 & \lambda^{-1}dz \\ 
      v(\lambda)\frac{n^2 z^{n-2}}{(z^n-1)^2}dz & 0 
    \end{pmatrix} \; , 
  \end{equation}
  where
  \begin{equation}
    \label{eq:vLemma}
      v(\lambda)=
      \frac{(n-2)^2w}{16n^2}(1-\lambda)^2 + 
  	\frac{1-n}{n^2}\lambda \; , \;\;\; 
      w\in \left[ \frac{-8 n}{(n-2)^2},\,0 
	\right) \; . 
  \end{equation}
  Let $w_0 \in \widetilde{\Sigma}$ be in the 
  fibre of $z_0=0 \in \Sigma$ and let $\Phi$ be the 
  solution of $d \Phi = \Phi \xi$, 
  $\Phi(w_0) = \Id$. Then there exists an initial 
  condition that unitarises the monodromy of 
  $\Phi$.
\end{theorem}
\begin{figure}[ht]
  \centering
  \includegraphics[scale=.8]
	{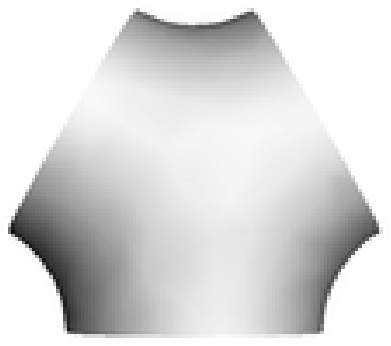}
  \includegraphics[scale=.8]
	{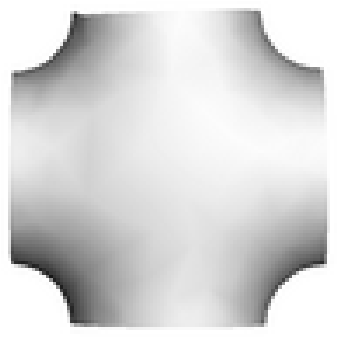}
  \includegraphics[scale=.8]
	{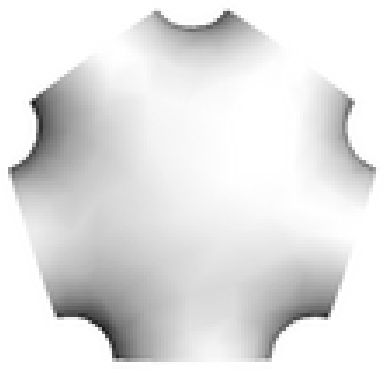}
  \includegraphics[scale=.8]
	{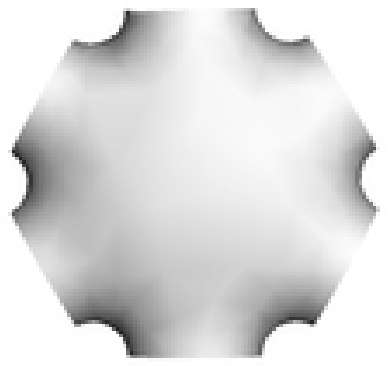}
  \includegraphics[scale=1.7]
	{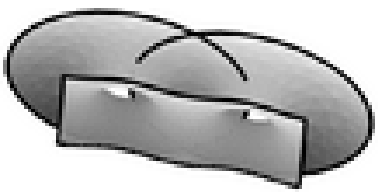}
\includegraphics[scale=1.1]
	{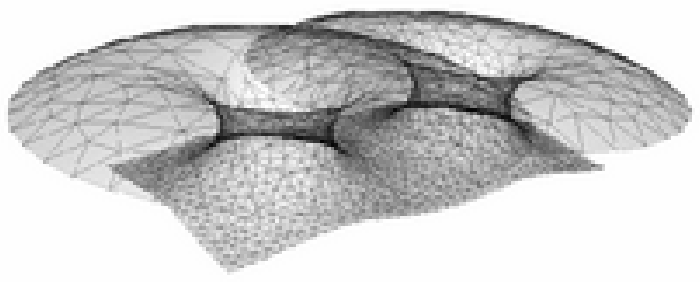}
\caption{{\sc{cmc}} surfaces with 
$3$-, $4$-, $5$- and 
$6$-fold symmetry (the {\sc{cmc}} immersions 
$f$ in \eqref{section5f} produced from 
Theorem \ref{lem:Section5first}) in the upper row.  
In the case of $3$-, $4$- and $6$-fold symmetry,
doubly periodic {\sc{cmc}} surfaces can be 
constructed 
by reflection in planes perpendicular to the 
plane of this page. The annulur ends of each
surface are nodoidal with equal weights and 
parallel same-directed axes.  
Two views of a portion of one of the 
doubly-periodic surfaces 
with $4$-fold symmetry are shown in the 
lower images.}
\label{fig:mattress}
\end{figure}

\begin{proof}
  Let $\alpha = \exp (\pi i/n)$, and define 
  the closed polygonal 
  loop $\gamma_0: [0,1] \to \Sigma$ as follows: 
\begin{equation*} 
  \gamma_0(t) = \left\{ 
  \begin{tabular}{llr}
   $4t\alpha^{-1}$ &for &$0\leq t\leq 1/4$, \\
   $(2-4 t)\alpha^{-1}+8t-2$ &for &$1/4\leq t\leq 1/2$, \\
   $6-8t+(4t-2)\alpha$ &for &$1/2\leq t\leq 3/4$, \\
   $(4-4t)\alpha$ &for &$3/4 \leq t \leq 1$. 
   \end{tabular}\right. 
\end{equation*}
 Then define the loops 
\begin{equation*} \gamma_j(t) = \alpha^{2j} \gamma_0(t) 
  \; , \quad j=1,2,...,n-1 \; . \end{equation*}  
  Let $M_j$ be the monodromy of $\Phi$ along 
  $\gamma_j$.  Then $M_0$, $M_1$, ..., $M_{n-1}$ 
  generate the monodromy 
  group of $\Phi$. Under the transformation 
  $\rho : z \to \alpha^2 z$ of $\Sigma$, we have 
  $\rho^\ast \xi = g^{-1}\,\xi \,g$, where 
  $g = \mathrm{diag}[\,\alpha^{-1},\,\alpha\,]$. 
  Because $\rho(z_0)=z_0$ and $\Phi(z_0) = \Id$, 
  we have $M_j = g^{-j} M_0 g^j$.  
  
  Changing variables to $\tilde z=1/z$ and 
  gauging ($\xi \mapsto \xi.g = 
  g^{-1}\,\xi\,g + g^{-1}\,dg$) by 
  $\tilde g = \mathrm{diag}[\,\tilde z^{-1},\, 
  \tilde z \,]$, we have 
\begin{equation*}
  \xi . \tilde g = 
  \begin{pmatrix} -\tilde z^{-1} & 
  -\lambda^{-1} \\ \frac{-v(\lambda) n^2 
  \tilde z^{n-2}}{(\tilde z^n-1)^2} 
   & \tilde z^{-1} \end{pmatrix} d\tilde z \; . 
\end{equation*} 
  Then one solution of $d\tilde \Phi = 
  \tilde \Phi \cdot (\xi . \tilde g)$ is 
  $\tilde \Phi = \exp \left( 
  \bigl( \begin{smallmatrix} -1 & 0 \\ 0 & 1 
  \end{smallmatrix} \bigr) \log \tilde z \right) 
     \tilde P(\tilde z,\lambda)$, 
  where $\tilde P (\tilde z,\lambda)$ is 
  well-defined and holomorphic with respect 
  to $\tilde z$ and is 
  nonsingular at $\tilde z=0$.  
  Furthermore, $\tilde P (\tilde z,\lambda)$ is 
  defined for all $\lambda \in \bbS^1$.  
  (This follows from a well-known 
  result in the theory of ordinary 
  differential equations, see \cite{KilKRS} 
  section 8.) It follows that the 
  monodromy of $\Phi$ along the loop 
  $\gamma_{n-1} ... \gamma_1 \gamma_0$ (here again 
  composition of these loops is from rightmost 
  first to leftmost last) encircling 
  $z=\infty$ is $(M_0) (g^{-1} M_0 g) ... 
  (g^{1-n} M_0 g^{n-1}) = \Id$.  Thus 
  $(M_0 g^{-1})^n = - \Id$.  Hence the 
  eigenvalues of $M_0 g^{-1}$ are constant 
  and are $n$-th roots of $-1$.  Since 
\begin{equation*} \Phi |_{\lambda = 1} = \begin{pmatrix} 
  d_z B & B \\ d_z D & D \end{pmatrix} \; , 
\end{equation*}
  with 
\begin{equation*} B = \alpha \sqrt[n]{z^n-1} \int_0^z 
         \left( \sqrt[n]{\zeta^n-1} \right)^{-2} 
  d\zeta\; , \;\;\; 
  D = \alpha^{-1} \sqrt[n]{z^n-1} \; , 
\end{equation*}
  we have that $M_0$ is upper-triangular 
  at $\lambda=1$ and the upper-left (resp. 
  lower-right) entry of its 
  diagonal is $\alpha^{-2}$ (resp. $\alpha^2$).  
  So the eigenvalues of 
  $M_0 g^{-1}$ are the same as the 
  eigenvalues of $g^{-1}$:
  \begin{equation} \label{gInverseEigvals}
     (\text{eigenvalues of $g^{-1}$}) = 
     (\text{eigenvalues of 
     $M_0 g^{-1}$}) = \alpha^{\pm 1} \; . 
  \end{equation}
  Now we determine the eigenvalues of $M_0$ 
  for general $\lambda$: For 
  $\hat g = \mathrm{diag}[\,\sqrt{z-1},\,
  \frac{1}{\sqrt{z-1}}\,]$ we have 
  \begin{equation*} \xi . \hat g = A \frac{dz}{z-1} + 
  O((z-1)^0) \; , \;\;\; A = \begin{pmatrix} \frac{1}{2} & 
             \lambda^{-1} \\ v(\lambda) & 
             \frac{-1}{2} \end{pmatrix} 
             \; , \end{equation*} in a neighborhood of $z=1$.  
  Applying Lemma 9.1 in \cite{KilKRS}, we have 
  that one solution of $d \hat \Phi = \hat \Phi 
  \cdot (\xi . \hat g )$ is 
  $\hat \Phi = \exp \left( A \,\log(z-1)\right) 
  \cdot \hat P (z,\lambda)$, where 
  $\hat P(z,\lambda)$ is holomorphic and 
  well-defined at $z=1$.  Furthermore, 
  $\hat P(z,\lambda)$ is 
  defined for any $\lambda \in \bbS^1$ at 
  which the difference of the eigenvalues of 
  $A$ is not an integer.  Hence 
  $\hat P(z,\lambda)$ is defined on 
  $\bbS^1$ minus a finite 
  set of points.  
  
  Hence one solution of 
  $d \check \Phi = \check \Phi \xi$ is 
  $\check \Phi = \hat \Phi \, \hat g^{-1}$.  
  Therefore any solution of $d \check \Phi = 
  \check \Phi \xi$ has monodromy along 
  $\gamma_0$ that is conjugate to 
  $- \exp (2 \pi i A)$.  In particular, the 
  eigenvalues of $M_0$ are $-\exp(\pm i\pi 
  \sqrt{1+4 \lambda^{-1} v(\lambda)})$, 
  and so 
  \begin{equation} \label{M0Eigvals}
     (\text{eigenvalues of $M_0$}) = 
     - \exp \left( \pm \tfrac{\pi i (n-2)}{n} 
   	\sqrt{1+\tfrac{w}{4} 
     \tfrac{(\lambda-1)^2}{\lambda}} \right) \; . 
  \end{equation}
  We now show that $M_0$ and $g$ can be 
  simultaneously unitarised at every 
  point in $\bbS^1$ where $\hat P(z,\lambda)$ 
  is defined: 
  We define the half-traces 
  $t_1=(1/2) \tr (g^{-1})$, 
  $t_2=(1/2) \tr (M_0 g^{-1})$ and 
  $t_3=(1/2) \tr (M_0)$.  
  Then, since $M_0\,g^{-1}\,(M_0 g^{-1})^{-1}=\Id$,   the condition for simultaneous unitarizability 
  \cite{Gol:top}, see also \cite{Bis}, of $M_0$ 
  and $g^{-1}$ and $(M_0 g^{-1})^{-1}$ is
  \begin{equation*} 
    1-t_1^2-t_2^2-t_3^2+2 t_1t_2t_3 \geq 0 \; . 
  \end{equation*}
  By Equations \eqref{gInverseEigvals} and 
  \eqref{M0Eigvals}, 
  this condition holds for all 
  $\lambda \in \bbS^1$ (where $\hat P(z,\lambda)$ 
  is defined) if and only if 
  \begin{equation*} 
    - \cos \left( \tfrac{\pi (n-2)}{n} 
      \sqrt{1+\tfrac{w}{4} 
       \tfrac{(\lambda-1)^2}{\lambda}} \right) \in 
       \left[ \cos \left( \tfrac{2 \pi}{n} \right),
  	\, 1 \right] \; , 
  \end{equation*}
  and this in turn holds if and only if 
  $w \in [\frac{-8n}{(n-2)^2},0]$, as in 
  Equation \eqref{eq:vLemma}.  
  
  It follows that the full monodromy group can 
  be unitarized at all but a 
  finite number of points in $\bbS^1$.  
  
  Note that if $M_0$ and $g^{-1}$ commute for 
  all $\lambda \in \bbS^1$, then $M_0$ must be 
  diagonal, and hence $M_0(\lambda) \in \SU$ 
  for all $\lambda \in \bbS^1$. 
  In this case, Lemma \ref{lem:Section5first} 
  is then clearly true, 
  so without loss of generality we may assume 
  that $[M_0,g^{-1}] \neq 0$.  
  Thus we can apply the gluing 
  theorem~\cite{Sch:tri} (see also \cite{KilKRS}) 
  to conclude there exists an initial condition  
  $h$ such that the monodromy group 
  of $h \Phi$ is unitary.  
\end{proof}

\begin{example} 
Now let $n$, $\calP$, $\Sigma$, $\xi$ and $w$ 
be as in Theorem~\ref{lem:Section5first}. 
Let $\calD=\{z \in \bbC \, | \, |z| \leq 1 \}$ 
be the closed unit disk in $\bbC$. 
Let $h = h(\lambda)$ be the unitariser of the 
monodromy of the solution $\Phi$ of 
$d\,\Phi = \Phi\,\xi$ given by 
Theorem~\ref{lem:Section5first}. Let 
\begin{equation}\label{section5f} 
      f:\calD \setminus \calP \to \bbR^3 
\end{equation} 
be the {\sc{cmc}} immersion generated by the data 
$(\,\Sigma,\,\xi,\,h,\,0\,)$. By the gluing theorem 
\cite{Sch:tri}, the immersion $f$ via the 
Sym-Bobenko formula \cite{Bob:tor}, 
in \eqref{section5f} is defined when using
$r$-Iwasawa splitting \cite{McI} for $r<1$ and $r$ 
sufficiently close to $1$. Then, up to a 
rigid motion and homothety of $\bbR^3$, 
we find numerically that $f$ has the following 
properties (see Figure~\ref{fig:mattress}): 
  \begin{itemize}
  \item the image of $f$ has order $n$ 
	dihedral symmetry, 
  \item the boundary of the image of $f$ 
	consists of $n$ complete planar 
        geodesics that are congruent to each 
	other, each lieing in a different plane 
  \begin{equation*} 
	\left\{ (x_1,x_2,x_3) \in \bbR^3 \, 
	\left| \, \cos \left( \tfrac{2\pi j}{n} 
          \right) x_1 + 
         \sin \left(\tfrac{2\pi j}{n} \right) 
	x_2 = 1 \right. \right\} 
  \end{equation*}
         for $j=0,1,...,n-1$, 
\item $f$ has $n$ ends at the punctures 
	in $\calP$, and the image of each end 
        is asymptotic to a $(\pi(n-2)/n)$-angle 
	arc of a Delaunay nodoid, 
\item the axes of the asymptotically 
	Delaunay ends are all vertical (i.e. 
        parallel to the line 
	$\{ (0,0,x_3) \in \bbR^3 \}$) and the 
        third coordinate $x_3$ of $f$ satisfies 
	$\lim_{z \in \calD, z \to p} x_3 = 
        +\infty$ for all $p \in \calP$.  
\item If $n\in\{3,\,4,\,6\}$, the complete 
	surface built by reflection across boundary 
        planar geodesics is doubly periodic; 
	in particular, it is invariant with 
	respect to two independent translations 
	of $\bbR^3$ parallel to the plane 
	$\{ (x_1,x_2,0) \in \bbR^3 \}$.  
  \end{itemize}
\end{example}

\begin{remark}
  In the cases $n=3,\,4,\,6$, the image 
  $f(\calD \setminus \calP)$ 
  can be repeatedly reflected to produce 
  a doubly-periodic surface with closed ends.  
  By the asymptotics theorem~\cite{Sch:tri}, 
  the annular ends are asymptotically Delaunay 
  with negative weight $w$.
\end{remark}

\section{Open problems}

\begin{enumerate}
\item Can one prove that the surfaces in 
Theorems \ref{thm:puncture1} and 
\ref{thm:doublypunctured} 
are complete and properly immersed?  
Could one further prove the 
asymptotic behavior of their ends?  
In particular, are the ends of the examples in 
Theorem \ref{thm:puncture1} asymptotic 
to ends of $2 n$-legged Smyth surfaces?  
\item By techniques like those used here, 
can one prove existence of 
a {\sc{cmc}} surface with finite topology 
and asymptotically Delaunay ends and 
positive genus?  
\end{enumerate}

%
%
\end{document}